\documentclass{amsart}
\usepackage[all]{xy}
\usepackage{color}

\newtheorem{theorem}{Theorem}[section]
\newtheorem{lemma}[theorem]{Lemma}
\newtheorem{proposition}[theorem]{Proposition}
\newtheorem{corollary}[theorem]{Corollary}

\newtheorem{definition}[theorem]{Definition}

\theoremstyle{remark}
\newtheorem{notation}[theorem]{Notation}

\theoremstyle{remark}
\newtheorem{remark}[theorem]{Remark}

\newcommand{\ko}{\: , \;}
\newcommand{\ca}{{\mathcal A}}

\newcommand{\cm}{{\mathcal M}}
\newcommand{\cb}{{\mathcal B}}
\newcommand{\cc}{{\mathcal C}}
\newcommand{\cq}{{\mathcal Q}}
\newcommand{\ch}{{\mathcal H}}
\newcommand{\cd}{{\mathcal D}}
\newcommand{\ci}{{\mathcal I}}
\newcommand{\cp}{{\mathcal P}}

\newcommand{\dgcat}{\mathsf{dgcat}}
\newcommand{\dgcatp}{\mathsf{dgcat}_{\geq 0}}

\newcommand{\internalcomment}[1]{}

\begin{document}

\title[DG versus Simplicial categories]{Differential graded versus Simplicial categories}
\author{Gon{\c c}alo Tabuada}
\address{Departamento de Matematica, FCT-UNL, Quinta da Torre, 2829-516 Caparica,~Portugal}


\email{
\begin{minipage}[t]{5cm}
tabuada@fct.unl.pt
\end{minipage}
}

\begin{abstract}
We construct a zig-zag of Quillen adjunctions between the homotopy theories of differential graded and simplicial categories. In an intermediate step we generalize Shipley-Schwede's work \cite{SS} on connective DG algebras by extending the Dold-Kan correspondence to a Quillen equivalence between categories enriched over positive graded chain complexes and simplicial $k$-modules. As an application we obtain a conceptual explanation of Simpson's homotopy fiber construction \cite{Simpson}.
\end{abstract}

\maketitle

\tableofcontents

\section{Introduction}
A differential graded (=dg) category is a category enriched in the category of complexes of modules over some commutative base ring
$k$. Dg categories provide a framework for `homological
geometry' and for `non commutative algebraic geometry' in the sense of
Drinfeld and Kontsevich \cite{Drinfeld} \cite{Chitalk} \cite{Kontsevich} \cite{IHP} \cite{finMotiv}.
In \cite{Tab} the homotopy theory of dg categories was constructed. This theory was allowed several developments such as: the creation by To{\"e}n of a derived Morita theory~\cite{Toen}; the construction of a category of non commutative motives~\cite{Tab}; the first conceptual characterization of Quillen-Waldhausen's $K$-theory~\cite{Tab}$\ldots$.

On the other hand a simplicial category is a category enriched over the category
of simplicial sets. Simplicial categories (and their close cousins: quasi-categories) provide a framework for `homotopy theories' and for `higher category theory' in the sense of Joyal, Lurie, Rezk, To{\"e}n $\ldots$ \cite{Joyal}\cite{Lurie}\cite{Rezk}\cite{Toen1}.
In \cite{Bergner} Bergner constructed a homotopy theory of simplicial categories by fixing an error in a previous version of \cite{Dwyer}. This theory can be considered as one of the four Quillen models for the theory of $(\infty, 1)$-categories, see \cite{survey} for a survey. 

We observe that the homotopy theories of differential graded and simplicial categories are formally similar and so a `bridge' between the two should be developed. In this paper we establish the first connexion between these theories by constructing a zig-zag of Quillen adjunctions relating the two:

In first place, we construct a Quillen model structure on positive graded dg categories by `truncating' the model structure of \cite{Tab}, see theorem~\ref{homotopypos}.

Secondly we generalize Shipley-Schwede's work \cite{SS} on connective DG algebras by extending the Dold-Kan correspondence to a Quillen equivalence between categories enriched over positive graded chain complexes and simplicial $k$-modules, see theorem~\ref{equivalenciafinal}.

Finally we extend the $k$-linearization functor to a Quillen adjunction between simplicial categories and simplicial $k$-linear categories.

As an application, the zig-zag of Quillen adjunctions obtained allow us to give a conceptual explanation of Simpson's homotopy fiber construction~\cite{Simpson} used in his nonabelian mixed Hodge theory.

\section{Acknowledgments}
I am deeply grateful to Gustavo Granja for several useful discussions and for his kindness. 

\section{Preliminaries}
In what follows, $k$ will denote a commutative ring with unit. The
tensor product $\otimes$ will denote the tensor product over $k$. Let
$\mathsf{Ch}$ denote the category of complexes over $k$ and $\mathsf{Ch}_{\geq 0}$ the full subcategory of positive graded complexes. Throughout this article we consider homological notation (the differential decreases the degree).  

Observe that $\mathsf{Ch}_{\geq 0}$ is a full symmetric monoidal subcategory of $\mathsf{Ch}$ and that the inclusion
$$ \mathsf{Ch}_{\geq 0} \hookrightarrow \mathsf{Ch}$$
commutes with limits and colimits. 

We denote by $\underline{\mathsf{Ch}_{\geq 0}}(-,-)$ the internal $\mbox{Hom}$-functor in $\mathsf{Ch}_{\geq 0}$ with respect to $\otimes$.

By a {\em dg category}, resp. {\em positive graded dg category}, we mean a category enriched over the symmetric monoidal category $\mathsf{Ch}$, resp. $\mathsf{Ch}_{\geq 0}$, see \cite{Drinfeld}~\cite{Keller94}~\cite{ICM}~\cite{Tab}. We denote by $\dgcat$, resp. $\dgcatp$, the category of small dg categories, resp. small positive graded dg categories. 

Notice that $\dgcatp$ is a full subcategory of $\dgcat$ and the inclusion
$$ \dgcatp \hookrightarrow \dgcat$$
commutes with limits and colimits.

Let $s\mathbf{Set}$ be the symmetric monoidal category of simplicial sets and $s\mathbf{Mod}$ the category of simplicial $k$-modules. We denote by $\wedge$ the levelwise tensor product of simplicial $k$-modules. The category $(s\mathbf{Mod}, -\wedge-)$ is a closed symmetric monoidal category.
We denote by $\underline{s\mathbf{Mod}}(-,-)$ its internal $\mbox{Hom}$-functor.

By a {\it simplicial category}, resp. {\em simplicial $k$-linear category}, we mean a category enriched over $s\mathbf{Set}$, resp. $s\mathbf{Mod}$, see \cite{Bergner}.

We denote by $s\mathbf{Set}$-$\mbox{Cat}$, resp. $s\mathbf{Mod}$-$\mbox{Cat}$, the category of small simplicial categories, resp. simplicial $k$-linear categories.

Let $(\cc, -\otimes-, \mathbb{I}_{\cc})$ and $(\cd, -\wedge-, \mathbb{I}_{\cd})$ be two symmetric monoidal categories. A {\em lax monoidal functor} is a functor $F:\cc \rightarrow \cd$ equipped with:
\begin{itemize}
\item[-] a morphism $\eta: \mathbb{I}_{\cd} \rightarrow F(\mathbb{I}_{\cc})$ and 
\item[-] natural morphisms
$$ \psi_{X,Y}: F(X) \wedge F(Y) \rightarrow F(X \otimes Y),\,\,\,\, X,Y \in \cc$$
which are coherently associative and unital (see diagrams $6.27$ and $6.28$ in \cite{Borceaux}).
\end{itemize}
A lax monoidal functor is {\em strong monoidal} if the morphisms $\eta$ and $\psi_{X,Y}$ are isomorphisms.

Throughout this article the adjunctions are displayed vertically with the left, resp. right, adjoint on the left side, resp. right side.
\section{Homotopy theory of positive graded DG categories}

In this section we will construct a Quillen model structure on $\mathsf{dgcat}_{\geq 0}$. For this we will adapt to our situation the Quillen model structure on $\mathsf{dgcat}$ constructed in chapter $1$ of \cite{Tab}.
\begin{remark}
Chapter $1$ of \cite{Tab} (and the whole thesis) is written using cohomological notation. Throughout this article we are always using homological notation.
\end{remark}
We now define the weak equivalences in $\dgcatp$.
\begin{definition}\label{quasieq}
A dg functor $F:\ca \rightarrow \cb$ in $\dgcatp$ is a {\em quasi-equivalence} if:
\item[(i)] $F(x,y): \ca(x,y) \rightarrow \cb(x,y)$ is a quasi-isomorphism in $\mathsf{Ch}_{\geq 0}$ for all objects $x,y \in \ca$ and
\item[(ii)] The induced functor $\mathsf{H}_0(F): \mathsf{H}_0(\ca) \rightarrow \mathsf{H}_0(\cb)$ is essentially surjective.
\end{definition}

\begin{notation}
We denote by $\cq_{qe}$ the class of quasi-equivalences in $\dgcatp$.
\end{notation}

\begin{remark}\label{quasi}
Notice that the class $\cq_{qe}$ consist exactly of those quasi-equivalences in $\dgcat$, see \cite[1.6]{Tab}, which belong to $\dgcatp$.
\end{remark}

In order to build a Quillen model structure on $\dgcatp$ we consider the generating (trivial) cofibrations in $\dgcat$ which belong to $\dgcatp$ and introduce a new generating cofibration. Let us now recall these constructions, see section $1.3$ in \cite{Tab}.

\begin{definition}\label{legos}
Following Drinfeld \cite[3.7.1]{Drinfeld} we define $\mathcal{K}$ to be the dg category that has two objects $1$, $2$ and whose morphisms are generated by $f \in \mathcal{K}(1,2)_0$,
$g \in \mathcal{K}(2,1)_0$,
$r_1 \in \mathcal{K}(1,1)_{1}$,
$r_2 \in \mathcal{K}(2,2)_{1}$ and $r_{12}
\in \mathcal{K}(1,2)_{2}$ subject to the relations $d(f)=d(g)=0$,
$d(r_1)=gf-\mathbf{1}_1$, $d(r_2) =fg-\mathbf{1}_2$ and $d(r_{12})=fr_1 - r_2f$.
$$\xymatrix{
    1 \ar@(ul,dl)[]_{r_1} \ar@/^/[r]^f \ar@/^0.8cm/[r]^{r_{12}} &
    2 \ar@(ur,dr)[]^{r_2} \ar@/^/[l]^g }
$$
Let $\underline{k}$ be the dg category with one object $3$, such that $\underline{k}(3,3)=k$. Let $F$ be the dg functor from $\underline{k}$ to $\mathcal{K}$ that sends $3$ to $1$.
Let $\mathcal{B}$ be the dg category with two objects $4$ and $5$ such that
$
\mathcal{B}(4,4)=k \ko
\mathcal{B}(5,5)=k \ko
\mathcal{B}(4,5)=0  
$
and $\mathcal{B}(5,4)=0$.
Let $n \geq 1$, $S^{n-1}$ the complex $k[n-1]$
and let $D^n$ be the mapping cone on the identity of $S^{n-1}$.
We denote by $\mathcal{P}(n)$ the dg category with two objects $6$ and $7$ such that
$
\mathcal{P}(n)(6,6)=k \ko
\mathcal{P}(n)(7,7)=k \ko
\mathcal{P}(n)(7,6)=0 \ko
\mathcal{P}(n)(6,7)=D^n
$ and whose composition given by multiplication.
Let $R(n)$ be the dg functor from $\mathcal{B}$ to
$\mathcal{P}(n)$ that sends $4$ to $6$ and $5$ to $7$.
Let $\mathcal{C}(n)$ be the dg category with two objects $8$ et $9$ such that
$
\mathcal{C}(n)(8,8)=k \ko
\mathcal{C}(n)(9,9)=k \ko
\mathcal{C}(n)(9,8)=0  \ko
\mathcal{C}(n)(8,9)=S^{n-1}
$ and whose composition given by multiplication.
Let $S(n)$ be the dg functor from $\mathcal{C}(n)$ to
$\mathcal{P}(n)$ that sends $8$ to $6$, $9$ to $7$ and $S^{n-1}$ to $D^n$
by the identity on $k$ in degree $n-1$.
Let $Q$ be the dg functor from the empty dg category $\emptyset$, which is the initial object in $\dgcatp$,
to $\underline{k}$. Finally let $N$ be the dg functor from $\cb$ to $\cc(1)$ that sends $4$ to $8$ and $5$ to $9$.
\end{definition}

Let us now recall the following standard recognition theorem:
\begin{theorem}{\cite[2.1.19]{Hovey}}\label{recognition}
Let $\cm$ be a complete and cocomplete category, $W$ a class of maps in $\cm$ and $I$ and $J$ sets of maps in $\cm$ such that:
\begin{itemize}
\item[1)] The class $W$ satisfies the two out of three axiom and is stable under retracts.
\item[2)] The domains of the elements of $I$ are small relative to $I$-$\mbox{cell}$.
\item[3)] The domains of the elements of $J$ are small relative to $J$-$\mbox{cell}$.
\item[4)] $J -\mbox{cell} \subseteq W\cap I-\mbox{cof}$.
\item[5)] $I-\mbox{inj}\subseteq W \cap J-\mbox{inj}$.
\item[6)] $W\cap I-\mbox{cof}\subseteq J-\mbox{cof}$ or
$W\cap J-\mbox{inj} \subseteq I-\mbox{inj}$.
\end{itemize}
Then there is a cofibrantly generated model category structure on $\cm$ in which $W$ is the class of weak equivalences, $I$ is a set of generating cofibrations, and $J$ is a set of generating trivial cofibrations.
\end{theorem}

\begin{theorem}\label{homotopypos}
If we let $\cm$ be the category $\dgcatp$, $W$ be the class $\cq_{qe}$, $J$ be the set of dg functors $F$ and $R(n), n \geq 1$, and $I$ the set of dg functors $Q$, $N$ and $S(n), n \geq 1$, then the conditions of the recognition theorem \ref{recognition} are satisfied. Thus, the category $\dgcatp$ admits a cofibrantly generated Quillen model structure whose weak equivalences are the quasi-equivalences.
\end{theorem}

\subsection{Proof of Theorem \ref{homotopypos}}
We start by observing that the category $\dgcatp$ is complete and cocomplete and that the class $\cq_{qe}$ satisfies the two out of three axiom and that it is stable under retracts. We observe also that the domains and codomains of the morphisms in $I$ and $J$ are small in the category $\dgcatp$. This implies that the first three conditions of the recognition theorem \ref{recognition} are verified.

\begin{lemma}\label{cell}
$J-\mbox{cell} \subseteq \cq_{qe}$.
\end{lemma}

\begin{proof}
Since the inclusion
$$ \dgcatp \hookrightarrow \dgcat$$
preserves colimits and the class $\cq_{qe}$ consist exactly of those quasi-equivalences in $\dgcat$ which belong to $\dgcatp$, the proof follows from lemma $1.10$ in \cite{Tab}.
\end{proof}
We now prove that $J-\mbox{inj} \cap \cq_{qe} = I-\mbox{inj}$. For this we introduce the following auxiliary class of dg functors:
\begin{definition}
Let $\mathbf{Surj}_{\geq 0}$ be the class of dg functors $G:\ch \rightarrow \ci$ in $\dgcatp$ such that:
\begin{itemize}
\item[-] $G(x,y): \ch(x,y) \rightarrow \ci(Gx,Gy)$ is a surjective quasi-isomorphism for all objects $x,y \in \ch$ and
\item[-] $G$ induces a surjective map on objects.
\end{itemize}
\end{definition}
\begin{remark}\label{dgcatpos}
Notice that the class $\mathbf{Surj}_{\geq 0}$ consist exactly of those dg functors in $\mathbf{Surj}$, see section $1.3.1$ in \cite{Tab}, which belong to $\dgcatp$.
\end{remark}

\begin{lemma}\label{classic}
$I-\mbox{inj} = \mathbf{Surj}_{\geq 0}$.
\end{lemma}

\begin{proof}
We prove first the inclusion $\supseteq$. Let $G: \ch \rightarrow \ci$ be a dg functor in $\mathbf{Surj}_{\geq 0}$. By remark~\ref{dgcatpos}, $G$ belongs to $\mathbf{Surj}$ and so lemma $1.11$ in \cite{Tab} implies that $G$ has the right lifting property with respect to the dg functors $Q$ and $S(n), n \geq 1$. Since the morphism of complexes
$$ G(x,y): \ch(x,y) \rightarrow \ci(Gx,Gy),\,\,\,\, x,y \in \ch$$
is surjective on the degree zero component, the dg functor $G$ also has the right lifting property with respect to $N$. This proves the inclusion $\supseteq$.

We now prove the inclusion $\subseteq$. Let $R: \cc \rightarrow \cd$ be a dg functor in $I-\mbox{inj}$. Lemma $1.11$ in \cite{Tab} implies that:
\begin{itemize}
\item[-] $R$ induces a surjective map on objects and
\item[-] for all objects $x,y \in \cc$:
\begin{itemize}
 \item[-] $R(x,y): \cc(x,y) \rightarrow \cd(Rx,Ry)$ is a surjective quasi-isomorphism for $n \geq 1$ and 
 \item[-] $\mathsf{H}_0 R(x,y): \mathsf{H}_0 \cc(x,y) \rightarrow \mathsf{H}_0 \cd(Rx,Ry)$ is an injective map.
\end{itemize} 
\end{itemize}
Since $R$ belongs to $I-\mbox{inj}$ it has the right lifting property with respect to $N$ and so the morphism of complexes $R(x,y)$ is also surjective on the degree zero component. This clearly implies that $R$ belongs to $\mathbf{Surj}_{\geq 0}$ and proves the inclusion $\subseteq$.
\end{proof}
We now consider the following `diagram chasing' lemma:
\begin{lemma}\label{chave}
Let $f:M_{\bullet} \rightarrow N_{\bullet}$ be a morphism in $\mathsf{Ch}_{\geq 0}$ such that:
\begin{itemize}
\item[-] $f_n: M_n \rightarrow N_n$ is surjective map for $n \geq 1$ and
\item[-] $\mathsf{H}_n(M_{\bullet}) \rightarrow \mathsf{H}_n(N_{\bullet})$ is an isomorphism for $n \geq 0$.
\end{itemize}
Then $f_0: M_0 \rightarrow N_0$ is also a surjective map.
\end{lemma}
\begin{proof}
It's a simple diagram chasing argument.
\end{proof}

\begin{lemma}\label{surjective}
$J-\mbox{inj} \cap \cq_{qe} = \mathbf{Surj}_{\geq 0}$.
\end{lemma}

\begin{proof}
The inclusion $\supseteq$ follows from remark~\ref{dgcatpos} and from the inclusion $\supseteq$ in lemma $1.12$ of \cite{Tab}. We now prove the inclusion $\subseteq$. Let $R:\cc \rightarrow \cd$ be a dg functor in $J-\mbox{inj}\cap \cq_{qe}$. Since $R$ belongs to $\cq_{qe}$ and it has the right lifting property with respect to the dg functors $R(n), n \geq 1$ the morphism of complexes
$$ R(x,y): \cc(x,y) \rightarrow \cd(Rx,Ry),\,\,\,\, x, y \in \cc$$
satisfies the conditions of lemma~\ref{chave} and so $R(x,y)$ is a surjective quasi-isomorphism. Finally the fact that $R$ induces a surjective map on objects follows from lemma $1.12$ in \cite{Tab}.

This proves the lemma.
\end{proof}

\begin{lemma}\label{implicit}
$J-\mbox{cell} \subseteq I-\mbox{cof}$.
\end{lemma}

\begin{proof}
Observe that the morphisms in $J-\mbox{cell}$ have the left lifting property with respect to the class $J-\mbox{inj}$. By lemmas ~\ref{classic} and \ref{surjective} $I-\mbox{inj}=J-\mbox{inj} \cap \cq_{qe}$ and so the morphisms in $J-\mbox{cell}$ have also the left lifting property with respect to the class $I-\mbox{inj}$, i.e. $J-\mbox{cell} \subseteq I-\mbox{cof}$.
\end{proof}

We have shown that $J-\mbox{cell} \subseteq \cq_{qe} \cap I-\mbox{cof} $ (lemmas~\ref{cell} and \ref{implicit}) and that $I-\mbox{inj} = J-\mbox{inj} \cap \cq_{qe}$ (lemmas~\ref{classic} and \ref{surjective}). This implies that the last three conditions of the recognition theorem \ref{recognition} are satisfied. This finishes the proof of theorem~\ref{homotopypos}.

\begin{remark}\label{fibrant}
Since every object in $\dgcat$ is fibrant, see remark $1.14$ in \cite{Tab}, and the set $J$ of generating trivial cofibrations in $\dgcatp$ is a subset of the generating trivial cofibrations in $\dgcat$ we conclude that every object in $\dgcatp$ is also fibrant.
\end{remark}

\subsection{The truncation functor}
In this subsection we construct a functorial path object in the Quillen model category $\dgcatp$.

Consider the following adjunction:
$$
\xymatrix{
\mathsf{Ch} \ar@<1ex>[d]^{\tau_{\geq 0}} \\
*+<1pc>{\mathsf{Ch}_{\geq 0}} \ar@{^{(}->}@<1ex>[u]^i ,
}
$$
where $\tau_{\geq 0}$ denotes the `intelligent' truncation functor: to a complex
$$ M_{\bullet}: \,\,\,\,\, \cdots \leftarrow M_{-2} \stackrel{d_{-1}}{\leftarrow} M_{-1} \stackrel{d_0}{\leftarrow} M_0 \stackrel{d_1}{\leftarrow} M_1 \leftarrow \cdots $$
it associates the complex
$$ \tau_{\geq 0}M_{\bullet}: \,\,\,\,\, \cdots \leftarrow 0 \leftarrow 0 \leftarrow Ker (d_0) \stackrel{d_1}{\leftarrow} M_1 \leftarrow \cdots\,.$$

The truncation functor $\tau_{\geq 0}$ is a lax monoidal functor. In particular we have natural morphisms
$$ \tau_{\geq 0} M_{\bullet} \otimes \tau_{\geq 0} N_{\bullet} \longrightarrow \tau_{\geq 0} (M_{\bullet} \otimes N_{\bullet}),\,\,\,\,M_{\bullet}, N_{\bullet} \in \mathsf{Ch} $$
which satisfy the associativity conditions. Observe that the truncation functor $\tau_{\geq 0}$ preserve the unit
$$ \cdots \leftarrow 0 \leftarrow k \leftarrow 0 \leftarrow \cdots$$
of both symmetric monoidal structures.

\begin{definition}
Let $\ca$ be a small dg category. The {\em truncation} $\tau_{\geq 0}\ca$ of  $\ca$ is the positive graded dg category with the same objects as $\ca$ and whose complexes of morphisms are defined as
$$ \tau_{\geq 0} \ca(x,y):= \tau_{\geq 0}\ca(x,y),\,\,\,\, x,y \in \ca \, .$$
For $x, y$ and $z$ objects in $\tau_{\geq 0} \ca$ the composition is defined as
$$\tau_{\geq 0} \ca(x, y) \otimes \tau_{\geq 0}\ca(y, z) \longrightarrow \tau_{\geq 0}(\ca(x, y)\otimes \ca(y, z)) \stackrel{\tau_{\geq 0}(c)}{\longrightarrow} \tau_{\geq 0}\ca(x, z),$$
where $c$ denotes the composition operation in $\ca$. The units in $\tau_{\geq 0}\ca$ are the same as those of $\ca$.
\end{definition}
Observe that we have a natural adjunction
$$
\xymatrix{
\dgcat \ar@<1ex>[d]^{\tau_{\geq 0}} \\
\dgcatp \ar@{^{(}->}@<1ex>[u]^i\,.
}
$$
\begin{remark}\label{eqadj}
Notice that both functors $i$ and $\tau_{\geq 0}$ preserve quasi-equivalences.
\end{remark}
\begin{proposition}\label{adjQuillen}
The adjunction $(i, \tau_{\geq 0})$ is a Quillen adjunction.
\end{proposition}

\begin{proof}
Clearly, by remark~\ref{quasi} the functor $i$ preserves weak equivalences. We now show that it also preserves cofibrations. The Quillen model structure of theorem~\ref{homotopypos} is cofibrantly generated and so by proposition $11.2.1$ in \cite{Hirschhornn} the class of cofibrations equals the class of retracts of relative $I-\mbox{cell}$ complexes. Since the functor $i$ preserves colimits it is then enough to prove that it sends the generating cofibrations in $\dgcatp$ to cofibrations in $\dgcat$. This is clear, by definition, for the generating cofibrations $Q$ and $S(n), n \geq 1$. We now observe that $i(N)=N$ is also a cofibration in $\dgcat$. In fact $N$ can be obtained by the following push-out
$$\xymatrix{
*+<1pc>{\cc(0)} \ar@{>->}[d]_{S(0)} \ar[r]^P  \ar@{}[dr]|{\lrcorner} & *+<1pc>{\cb} \ar@{>->}[d]^N \\
\cp(0) \ar[r] & \cc(1)\,,
}
$$
where $S(0)$ is a generating cofibration in $\dgcat$, see section $1.3$ in \cite{Tab}, and $P$ sends $8$ to $4$ and $9$ to $5$. This proves the lemma.
\end{proof}
\begin{remark}
Recall from \cite[4.1]{Tab} the construction of a path object $P(\ca)$ for each dg category $\ca \in \dgcat$.
\end{remark}
\begin{lemma}\label{path}
Let $\ca$ be a positive graded dg category. Then $\tau_{\geq 0} P(\ca)$ is a path object of $\ca$ in $\dgcatp$.
\end{lemma}
\begin{proof}
Consider the diagonal dg functor
$$ \ca \stackrel{\Delta}{\longrightarrow} \ca \times \ca$$
in $\dgcat$. We have, as in \cite[4.1]{Tab}, a factorization 
$$ 
\xymatrix{
\ca \ar[rr]^{\Delta} \ar[dr]_I^{\sim} && \ca \times \ca \\
& P(\ca) \ar@{->>}[ur]_P \,,& 
}
$$
where $I$ is a quasi-equivalence and $P$ a fibration. By remark~\ref{eqadj} and lemma~\ref{adjQuillen} the functor $\tau_{\geq 0}$ preserves quasi-equivalences and fibrations. Since the functor $\tau_{\geq 0}$ also preserves limits we obtain the following factorization
$$
\xymatrix{
\ca \ar[rr]^{\Delta} \ar[dr]^{\sim}_{\tau_{\geq 0}(I)} & & \ca \times \ca \\
& \tau_{\geq 0} P(\ca) \ar@{->>}[ur]_{\tau_{\geq 0}(P)} & .
}
$$
This proves the lemma.
\end{proof}

\section{Extended Dold-Kan equivalence}
In this section we will first construct a Quillen model structure on $s\mathbf{Mod}$-$\mbox{Cat}$ and then show that it is Quillen equivalent to the model structure on $\dgcatp$ of theorem~\ref{homotopypos}.

Recall from \cite[III-2.3]{Jardine} the Dold-Kan equivalence between simplicial $k$-modules and positive graded complexes
$$
\xymatrix{
s\mathbf{Mod} \ar@<1ex>[d]^N \\
\mathsf{Ch}_{\geq 0} \ar@<1ex>[u]^{\Gamma} ,
}
$$
where $N$ is the normalization functor and $\Gamma$ its inverse. The normalization functor $N$ is a lax monoidal functor, see~\cite[2.3]{SS}, via the Eilenberg-MacLane shuffle map, see \cite[VIII-8.8]{Maclane}
$$ \nabla :NA \otimes NB \longrightarrow N(A\wedge B), \,\,\,\, A,B \in s\mathbf{Mod}\,.$$
Observe that the normalization functor $N$ preserves the unit of the two symmetric monoidal structures.

As it is shown in \cite[2.3]{SS} the lax monoidal structure on $N$, given by the shuffle map $\nabla$, induces a lax {\em co}monoidal structure on $\Gamma$:
$$ \widetilde{\psi}: \Gamma(M\otimes M') \longrightarrow \Gamma(M) \wedge \Gamma(M'),\,\,\,\, M, M' \in \mathsf{Ch}_{\geq 0}.$$

Now, let $I$ be a set. 
\begin{notation}
We denote by $\mathsf{Ch}_{\geq 0}^I$-$\mbox{Gr}$, resp. by $\mathsf{Ch}_{\geq 0}^I$-$\mbox{Cat}$, the category of $\mathsf{Ch}_{\geq 0}$-graphs with a fixed set of objects $I$, resp. the category of categories enriched over $\mathsf{Ch}_{\geq 0}$ which have a fixed set of objects $I$. The morphisms in $\mathsf{Ch}_{\geq 0}^I$-$\mbox{Gr}$ and $\mathsf{Ch}_{\geq 0}^I$-$\mbox{Cat}$ induce the identity map on the objects.
\end{notation}
We have a natural adjunction
$$
\xymatrix{
\mathsf{Ch}_{\geq 0}^I\text{-}\mbox{Cat} \ar@<1ex>[d]^U \\
\mathsf{Ch}_{\geq 0}^I \text{-}\mbox{Gr} \ar@<1ex>[u]^{T_I},
}
$$
where  $U$ is the forgetful functor and $T_I$ is defined as
\[ T_I(\ca)(x,y):= \left\{ \begin{array}{ll}
k \oplus \underset{x, x_1, \ldots, x_n, y}{\bigoplus} \ca(x,x_1)  \otimes \ldots \otimes \ca(x_n,y) & \mbox{if $x=y$} \\
\underset{x, x_1, \ldots, x_n,y}{\bigoplus} \ca(x,x_1) \otimes \ldots \otimes \ca(x_n,y) & \mbox{if $x\neq y$} \end{array} \right. \]
Composition is given by concatenation and the unit corresponds to $\mathbf{1} \in k$.

\begin{remark}\label{Qmodel}
\begin{itemize}
\item[-] Notice that the categories $\mathsf{Ch}_{\geq 0}^I$-$\mbox{Gr}$ and $\mathsf{Ch}_{\geq 0}^I$-$\mbox{Cat}$ admit standard Quillen model structures whose weak equivalences (resp. fibrations) are the morphisms $F:\ca \rightarrow \cb$ such that
$$ F(x,y):\ca(x,y) \longrightarrow \cb(x,y),\,\,\,\, x,y \in I$$
is a weak equivalence (resp. fibration) in $\mathsf{Ch}_{\geq 0}$. In fact the projective Quillen model structure on $\mathsf{Ch}_{\geq 0}$, see~\cite[III-2]{Jardine}, naturally induces a model structure on $\mathsf{Ch}_{\geq 0}^I\text{-}\mbox{Gr}$ which can be lifted along the functor $T_I$ using theorem $11.3.2$ in \cite{Hirschhornn}.
\item[-] If the set $I$ has a unique element, then the previous adjunction corresponds to the (Quillen) adjunction between connective dg algebras and positive graded complexes, see~\cite{Jardine1}.
\end{itemize}
\end{remark}

\begin{notation}
We denote by $s\mathbf{Mod}^I$-$\mbox{Gr}$, resp. by $s\mathbf{Mod}^I$-$\mbox{Cat}$, the category of $s\mathbf{Mod}$-graphs with a fixed set of objects $I$, resp. the category of categories enriched over $s\mathbf{Mod}$ which have a fixed set of objects $I$. The morphisms in $s\mathbf{Mod}^I$-$\mbox{Gr}$ and $s\mathbf{Mod}^I$-$\mbox{Cat}$ induce the identity map on the objects.
\end{notation}

In an analogous way we have an adjunction
$$
\xymatrix{
s\mathbf{Mod}^I\text{-}\mbox{Cat} \ar@<1ex>[d]^U \\
s\mathbf{Mod}^I\text{-}\mbox{Gr} \ar@<1ex>[u]^{T_I},
}
$$
where  $U$ is the forgetful functor and $T_I$ is defined as
\[ T_I(\cb)(x,y):= \left\{ \begin{array}{ll}
 k\Delta_0 \oplus \underset{x,x_1,\ldots,x_n,y}{\bigoplus} B(x,x_1)  \wedge \ldots \wedge B(x_n,y)& \mbox{if $x=y$} \\
 \underset{x, x_1, \ldots, x_n,y}{\bigoplus} B(x,x_1)  \wedge \ldots \wedge B(x_n,y) & \mbox{if $x\neq y$} \end{array} \right. \]
Composition is given by concatenation and the unit corresponds to $\mathbf{1} \in k\Delta_0$.
\begin{remark}
If the set $I$ has an unique element, then the previous adjunction corresponds to the classical adjunction between simplicial $k$-algebras and simplicial $k$-modules, see~\cite[II-5.2]{Jardine}.
\end{remark}
Clearly the Dold-Kan equivalence induces an equivalence of categories
$$
\xymatrix{
s\mathbf{Mod}^I\text{-}\mbox{Gr} \ar@<1ex>[d]^N \\
\mathsf{Ch}_{\geq 0}^I\text{-}\mbox{Gr} \ar@<1ex>[u]^{\Gamma}
}
$$
that we still denote by $N$ and $\Gamma$.

Since the functor $N: s\mathbf{Mod} \rightarrow \mathsf{Ch}_{\geq 0}$ is lax monoidal it induces, as in \cite[3.3]{SS}, a normalization functor
$$
\xymatrix{
s\mathbf{Mod}^I\text{-}\mbox{Cat} \ar@<1ex>[d]^{N_I} \\
\mathsf{Ch}_{\geq 0}^I\text{-}\mbox{Cat}\,.
}
$$
In fact, let $A \in s\mathbf{Mod}^I$-$\mbox{Cat}$ and $x, y$ and $z$ objects in $A$. Then $N_I(A)$ has the same objects as $A$, the complexes of morphisms are given by
$$ N_I(A)(x,y):=NA(x,y),\,\,\,\, x, y \in A$$
and the composition is defined by
$$ NA(x,y) \otimes NA(y,z) \stackrel{\nabla}{\longrightarrow} N(A(x,y)\wedge A(y,z)) \stackrel{N(c)}{\longrightarrow} NA(x,z)\,,$$
where $c$ denotes the composition operation in $A$. The units in $N_I(A)$ are induced by those of $A$ under the normalization functor $N$.

As it is shown in section $3.3$ of \cite{SS} the functor $N_I$ admits a left adjoint $L_I$.

Let $\ca \in  \mathsf{Ch}_{\geq 0}^I$-$\mbox{Cat}$. The value of the left adjoint $L_I$ on $\ca$ is defined as the coequalizer of two morphisms in $s\mathbf{Mod}^I$-$\mbox{Cat}$
$$ \xymatrix{ T_I \Gamma U T_IU(\ca) \ar@<1ex>[r]^{\psi_1} \ar@<-1ex>[r]_{\psi_2}  & T_I \Gamma U(\ca) \ar[r] & L_I(\ca).}$$
The morphism $\psi_1$ is obtained from the unit of the adjunction
$$ T_I U \ca \longrightarrow \ca$$
by applying the composite functor $T_I \Gamma U$; the morphism $\psi_2$ is the unique morphism in $s\mathbf{Mod}^I$-$\mbox{Cat}$ induced by the $s\mathbf{Mod}^I$-$\mbox{Gr}$ morphism
$$ \Gamma U T_I U(\ca) \longrightarrow U T_I \Gamma U(\ca)$$
whose value at $\Gamma U T_I U(\ca)(x,y),\,\,\,\, x, y \in I$ is
$$ 
\xymatrix{
\underset{x,x_1, \ldots, x_n,y}{\bigoplus} \Gamma(\ca(x,x_1)\otimes \ldots \otimes \ca(x_n,y)) \ar[d]^-{\tilde{\psi}} \\ \underset{x,x_1, \ldots, x_n,y}{\bigoplus} \Gamma \ca(x,x_1) \wedge \ldots \wedge \Gamma \ca(x_n,y),
}
$$
where $\tilde{\psi}$ is the lax {\em co}monoidal structure on $\Gamma$ induced by the lax monoidal structure on $N$, see section $3.3$ of \cite{SS}.

\subsection{Left adjoint}
Notice that the normalization functor $N_I: s\mathbf{Mod}^I$-$\mbox{Cat} \rightarrow \mathsf{Ch}_{\geq 0}^I$-$\mbox{Cat}$, of the previous subsection, can be naturally defined for every set $I$ and so it induces a `global' normalization functor
$$
\xymatrix{
s\mathbf{Mod}\text{-}\mbox{Cat} \ar@<1ex>[d]^N \\
\dgcatp .
}
$$
In this subsection we will construct the left adjoint of $N$.

Let $\ca \in \dgcatp$ and denote by $I$ its set of objects. Define $L(\ca)$ as the simplicial $k$-linear category $L_I(\ca)$. 

Now, let $F:\ca \rightarrow \ca'$ be a dg functor. We denote by $I'$ the set of objects of $\ca'$. The dg functor $F$ induces the following diagram in $s\mathbf{Mod}\text{-}\mbox{Cat}$:
$$
\xymatrix{
T_I\Gamma U T_I U(\ca) \ar@<1ex>[r]^{\psi_1} \ar@<-1ex>[r]_{\psi_2} \ar[d] & T_I \Gamma U(\ca) \ar[d] \ar[r] & L_I(\ca)=:L(\ca)  \\
T_{I'}\Gamma U T_{I'} U(\ca') \ar@<1ex>[r]^{\psi_1} \ar@<-1ex>[r]_{\psi_2} & T_{I'} \Gamma U(\ca') \ar[r] & L_{I'}(\ca')=:L(\ca')\,.
}
$$
Notice that the square whose horizontal arrows are $\psi_1$ (resp. $\psi_2$) is commutative. Since the inclusions
$$ s\mathbf{Mod}^I\text{-}\mbox{Cat} \hookrightarrow s\mathbf{Mod}\text{-}\mbox{Cat} \,\,\,\,\text{and} \,\,\,\, s\mathbf{Mod}^{I'}\text{-}\mbox{Cat} \hookrightarrow s\mathbf{Mod}\text{-}\mbox{Cat}$$
clearly preserve coequalizers the previous diagram in $s\mathbf{Mod}\text{-}\mbox{Cat}$ induces a simplicial $k$-linear functor
$$ L(F) : L(\ca) \longrightarrow L(\ca').$$
We have constructed a functor
$$ L: \dgcatp \longrightarrow s\mathbf{Mod}\text{-}\mbox{Cat}.$$

\begin{proposition}
The functor $L$ is left adjoint to $N$. 
\end{proposition}

\begin{proof}
Let $\ca \in \dgcatp$ and $B \in s\mathbf{Mod}\text{-}\mbox{Cat}$. Let us denote by $I$ the set of objects of $\ca$. We will construct two natural maps
$$\xymatrix{
 s\mathbf{Mod}\text{-}\mbox{Cat}(L(\ca),B) \ar@<1ex>[r]^{\phi} & \dgcatp (\ca, N(B)) \ar@<1ex>[l]^{\eta}
 }$$
and then show that they are inverse of each other.

Let $G: L(\ca) \rightarrow B$ be a simplicial $k$-linear functor. We denote by $B'$ the full subcategory of $B$ whose objects are those which belong to the image of $G$. We have a natural factorization
$$
\xymatrix{
L(\ca) \ar[rr]^{G} \ar[dr]_{G'} &  & B \\
 & B' \ar@{^{(}->}[ur]& .
}
$$
Now, let $\widetilde{B}$ be the simplicial $k$-linear category whose set of objects is
$$\mbox{obj}(\widetilde{B}):=\{ (a,b) |\,a \in L(\ca), b \in B' \, \mbox{and}\, G'(a)=b \}$$
and whose simplicial $k$-module of morphisms is defined as
$$ \widetilde{B}((a,b),(a',b')):= B'(b,b').$$
The composition is given by the composition in $B'$. Now consider the simplicial $k$-linear functor
$$ \widetilde{G} : L(\ca) \longrightarrow \widetilde{B}$$
which maps $a$ to $(a,G'(a))$ and the simplicial $k$-linear functor 
$$ P: \widetilde{B} \longrightarrow B'$$
which maps $(a,b)$ to $b$. 

The above constructions allow us to factor $G$ as the following composition
$$
\xymatrix{
L(\ca) \ar[r]^{\widetilde{G}} & \widetilde{B} \ar[r]^P & B' \ar@{^{(}->}[r] & B.
}
$$
Notice that $\widetilde{G}$ induces a bijection on objects and so it belongs to $s\mathbf{Mod}^I\text{-}\mbox{Cat}$. Finally define $\phi(G)$ as the following composition
$$
\xymatrix{
\phi(G):\ca \ar[r]^-{\widetilde{G}^{\sharp}} & N\widetilde{B} \ar[r]^{NP} &  NB' \ar@{^{(}->}[r] & NB,
}
$$
where $\widetilde{G}^{\sharp}$ denotes the morphism in $\mathsf{Ch}_{\geq 0}^I\text{-}\mbox{Cat}$ which corresponds to $\widetilde{G}$ under the adjunction $(L_I, N_I)$.

We now construct in a similar way the map $\eta$. Let $F:\ca \rightarrow NB$ be a dg functor and $(NB)'$ be the full subcategory of $NB$ whose objects are those which belong to the image of $F$. We have a natural factorization
$$
\xymatrix{
\ca \ar[rr]^{F} \ar[dr]_{F'} &  & NB \\
 & (NB)' \ar@{^{(}->}[ur]& .
}
$$
Now, let $\widetilde{NB}$ be the positive graded dg category whose set of objects is
$$\mbox{obj}(\widetilde{NB}):=\{ (a,b) |\, a \in \ca, b\in NB' \, \mbox{and}\,  F'(a)=b \}$$
and whose positive graded complex of morphisms is defined as
$$ \widetilde{NB}((a,b),(a',b')):= (NB)'(b,b').$$
The composition is given by the composition in $(NB)'$. Consider the dg functor
$$\widetilde{F}: \ca \longrightarrow \widetilde{NB}$$
which maps $a$ to $(a,F'(a))$ and the dg functor
$$ P: \widetilde{NB} \longrightarrow (NB)'$$
which maps $(a,b)$ to $b$.

The above constructions allow us to factor $F$ as the following composition
$$ 
\xymatrix{
\ca \ar[r]^{\widetilde{F}} & \widetilde{NB} \ar[r]^-P & (NB)' \ar@{^{(}->}[r] & NB.
}
$$
Notice that $\widetilde{F}$ induces a bijection on objects and so belongs to $\mathsf{Ch}_{\geq 0}^I\text{-}\mbox{Cat}$. Since the normalization functor $N$ preserve the set of objects, the above construction
$$ 
\xymatrix{
\widetilde{NB} \ar[r]^-P & (NB)' \ar@{^{(}->}[r] & NB
}
$$
can be naturally lifted to the category $s\mathbf{Mod}\text{-}\mbox{Cat}$. We have the folowing diagram
$$
\xymatrix{
s\mathbf{Mod}\text{-}\mbox{Cat} \ar[d]^N & \widetilde{B} \ar[r]^{\overline{P}} \ar@{|->}[d] & *+<1pc>{B'} \ar@{|->}[d] \ar@{^{(}->}[r] & B \ar@{|->}[d] \\
\dgcatp & \widetilde{NB} \ar[r]^P & (NB)' \ar@{^{(}->}[r] & NB.
}
$$
We can now define $\eta(F)$ as the following composition
$$ 
\xymatrix{
\eta(F): L(\ca) \ar[r]^-{\widetilde{F}^{\natural}} & \widetilde{B} \ar[r]^{\overline{P}} & B' \ar@{^{(}->}[r] & B, 
}
$$
where $\widetilde{F}^{\natural}$ denotes the morphism in $s\mathbf{Mod}^I\text{-}\mbox{Cat}$, which corresponds to $\widetilde{F}$ under the adjunction $(L_I, N_I)$.

The maps $\eta$ and $\phi$ are clearly inverse of each other and so the proposition is proven. 
\end{proof}

\subsection{Path object}
In this subsection we lift the Quillen model structure on $\dgcatp$, see theorem~\ref{homotopypos}, along the adjunction
$$
\xymatrix{
s\mathbf{Mod}\text{-}\mbox{Cat} \ar@<1ex>[d]^N \\
\dgcatp \ar@<1ex>[u]^L
}
$$
of the previous subsection.
For this we will use theorem $5.12$ and proposition $5.13$ of \cite{Tab}.

\begin{definition}\label{abmorph}
A simplicial $k$-linear functor $G: A \rightarrow B$ is:
\begin{itemize}
\item[-] a {\em weak equivalence} if $NG$ is a quasi-equivalence in $\dgcatp$.
\item[-] a {\em fibration} if $NG$ is a fibration in $\dgcatp$.
\item[-] a {\em cofibration} if it has the left lifting property with respect to all trivial fibrations in $s\mathbf{Mod}\text{-}\mbox{Cat}$.
\end{itemize}
\end{definition}

\begin{definition}
Let $A$ be a small simplicial $k$-linear category. The {\em homotopy category} $\pi_0(A)$ of $A$ is the category which has the same objects as $A$ and whose morphisms are defined as 
$$\pi_0(A)(x,y):= \pi_0(A(x,y)), \,\,\, x, y \in A\,.$$
\end{definition}

\begin{lemma}\label{key}
Let $x \stackrel{f}{\rightarrow} y$ be a $0$-simplex morphism in $A$. Then $\pi_0(f)$ is invertible in $\pi_0(A)$ iff $\mathsf{H}_0(Nf)$ is invertible in $\mathsf{H}_0(NA)$.
\end{lemma}
\begin{proof}
We start by observing that if we restrict ourselves to the $0$-simplex morphisms in $A$ and to the degree zero morphisms in $NA$ we have the same category. In fact the degree zero component of the shuffle map $\nabla$, used in the definition of $NA$, is the identity map, see \cite[VIII-8.8]{Maclane}.

Now suppose that $\pi_0(f)$ is invertible. Then there exists a $0$-simplex morphism $g:y \rightarrow x$ and $1$-simplex morphisms $h_1 \in A(x,x)$ and $h_2 \in A(y,y)$ such that $d_0(h_1)=\mathbf{1}_X, d_1(h_1)=gf, d_0(h_2)=\mathbf{1}_Y$ and $d_1(h_2)=fg$. Observe that the image of $h_1$, resp. $h_2$, by the normalization functor $N$ is a degree $1$ morphism in $NA(x,x)$, resp. in $NA(y,y)$, whose differential is $gf-\mathbf{1}_X$ (resp. $fg-\mathbf{1}_Y$). This implies that $\mathsf{H}_0(Nf)$ is also invertible in $\mathsf{H}_0(NA)$.

To prove the converse we consider an analogous argument.
\end{proof}

\begin{proposition}\label{eqfracas}
A simplicial $k$-linear functor $G:A \rightarrow B$ is a weak equivalence iff:
\begin{itemize}
\item[(1)] $G(x,y) :A(x,y) \rightarrow B(Gx,Gy)$ induces an isomorphism on $\pi_i$ for all $i \geq 0$ and for all objects $x,y \in A$, and
\item[(2)] $\pi_0(G): \pi_0(A) \rightarrow \pi_0(B)$ is essentially surjective.
\end{itemize}
\end{proposition}

\begin{proof}
We show that condition $(1)$, resp. condition $(2)$, is equivalent to condition $(i)$, resp. condition $(ii)$, of definition~\ref{quasieq}. By the Dold-Kan equivalence, we have the following commutative diagram
$$ \xymatrix{
\pi_iA(x,y) \ar[r]^G \ar[d]^-{\sim} & \pi_iB(Gx,Gy) \ar[d]^{\sim} \\
\mathsf{H}_iNA(x,y) \ar[r]_-{NG} & \mathsf{H}_iNB(Gx,Gy)
}
$$
where the vertical arrows are isomorphisms. This implies that condition $(1)$ is equivalent to condition $(i)$ of definition~\ref{quasieq}.

Concerning condition $(2)$, we start by supposing that $\pi_0(G)$ is essentially surjective. Consider the functor
$$ \mathsf{H}_0(NG): \mathsf{H}_0(NA) \rightarrow \mathsf{H}_0(NB)$$
and let $z$ be an object in $\mathsf{H}_0(NB)$. Since $\pi_0(B)$ and $\mathsf{H}_0(NB)$ have the same objects we can consider $z$ as an object in $\pi_0(B)$. By hypothesis, $\pi_0(G)$ is essentially surjective and so there exists an object $w \in \pi_0(A)$ and a $0$-simplex morphism
$$ Gw \stackrel{f}{\rightarrow} z$$
which becomes invertible in $\pi_0(B)$. Now lemma~\ref{key} implies that $Nf$ is invertible in $\mathsf{H}_0(NB)$ and so we conclude that the functor $\mathsf{H}_0(NG)$ is essentially surjective. This shows that condition $(2)$ implies condition $(ii)$ of definition~\ref{quasieq}. To prove the converse we consider an analogous argument.
\end{proof}

\begin{theorem}\label{newQuillen}
The category $s\mathbf{Mod}\text{-}\mbox{Cat}$ when endowed with the notions of weak equivalence, fibration and cofibration as in definition~\ref{abmorph}, becomes a cofibrantly generated Quillen model category and the adjunction $(L,N)$ becomes a Quillen adjunction.
\end{theorem}

The proof will consist on verifying the conditions of theorem $5.12$ and proposition $5.13$ in \cite{Tab}. Since the Quillen model structure on $\dgcatp$ is cofibrantly generated, see theorem~\ref{homotopypos}; every object in $\dgcatp$ is fibrant, see remark~\ref{fibrant}; and the functor $N$ clearly commutes with filtered colimits it is enough to show that:
\begin{itemize}
\item[-] for each simplicial $k$-linear category $A$, we have a factorization
$$
\xymatrix{
A \ar[rr]^{\Delta} \ar[dr]^{\sim}_{I_A} & & A \times A \\
& P(A) \ar@{->>}[ur]_{P_0 \times P_1}, & 
}
$$
with $I_A$ is a weak equivalence and $P_0 \times P_1$ is a fibration in $s\mathbf{Mod}\text{-}\mbox{Cat}$. 
\end{itemize}
For this we need a few lemmas. We start with the following definition.
\begin{definition}
Let us define $P(A)$ as the simplicial $k$-linear category whose objects are the $0$-simplex morphisms $f:x\rightarrow y$ in $A$ which become invertible in $\pi_0(A)$. We define the simplicial $k$-module of morphisms
$$ P(A)(x \stackrel{f}{\rightarrow}y, x' \stackrel{f'}{\rightarrow} y'),\,\,\,\, f, f' \in P(A)$$
as the homotopy pull-back in $s\mathbf{Mod}$ of the diagram
$$
\xymatrix{
& A(y,y') \ar[d]^{f^{\ast}} \\
A(x,x') \ar[r]_{f'_{\ast}} & A(x,y')\,,
}
$$
by which we mean the simplicial $k$-module
$$ A(x,x') \underset{A(x,y')}{\times} \underline{s\mathbf{Mod}}(k\Delta[1], A(x,y')) \underset{A(x,y')}{\times} A(y,y')\,.$$
We denote the simplexes in $A(x,x')$ and $A(y,y')$ {\em lateral morphisms} and the simplexes in $s\mathbf{Mod}(k\Delta[1], A(x,y'))$ {\em homotopies}.
The composition operation
$$ P(A)(f,f') \wedge P(A)(f',f'') \longrightarrow P(A)(f,f''), \,\,\,\, f, f', f'' \in P(A)$$
decomposes on:
\begin{itemize}
\item[-] a composition of lateral morphisms, which is induced by the composition on $A$ and
\item[-] a composition of homotopies, which is given by the map
$$ 
\xymatrix@!0 @R=3,5pc @C=1pc{
\underline{s\mathbf{Mod}}(k\Delta[1],A(x,y'))\wedge A(y',y'') \underset{\underline{s\mathbf{Mod}}(k\Delta[0], A(x,y''))}{\times} A(x,x') \wedge \underline{s\mathbf{Mod}}(k\Delta[1], A(x',x'')) \ar[d]^{\text{composition}} \\
\underline{s\mathbf{Mod}}(k\Delta[1] \underset{k\Delta[0]}{\bigoplus} k\Delta[1], A(x,y'')) \ar[d] \\
\underline{s\mathbf{Mod}}(k\Delta[1], A(x,y''))\,,
}
$$
where the last map is induced by the diagonal map in $k\Delta[1]$.
\end{itemize}
\end{definition}
\begin{remark}\label{nova}
Notice that a $0$-simplex morphism $\alpha:f \rightarrow f'$ in $P(A)$ is of the form $\alpha= (m_x, h, m_y)$, with $m_x: x\rightarrow x'$ and $m_y:y \rightarrow y'$ $0$-simplex morphisms in $A$ and $h$ is a $1$-simplex morphism in $A(x,y')$ such that $d_0(h)=m_yf$ and $d_1(h)=f'm_x$.
\end{remark}
We have a natural commutative diagram in $s\mathbf{Mod}\text{-}\mbox{Cat}$
\begin{equation}\label{numero}
\xymatrix{
A \ar[rr]^{\Delta} \ar[dr]_{I_A} && A \times A \\
& P(A) \ar[ur]_{P_0\times P_1}
}
\end{equation}
where $I_A$ is the simplicial $k$-linear functor that associates to an object $x \in A$ the $0$-simplex morphism $x \stackrel{Id}{\rightarrow} x$ and $P_0$, resp. $P_1$, is the simplicial $k$-linear functor that sends a morphism $x \stackrel{f}{\rightarrow} y$ in $P(A)$ to $x$, resp. $y$.

Notice that by applying the normalization functor $N$ to the above diagram and lemma~\ref{path} to the dg category $NA$ we obtain two factorizations
$$
\xymatrix{
NA \ar[rr]^{\Delta} \ar[dr] \ar@/_1pc/[ddr]_{\tau_{\geq 0}(I)} && NA \times NA \\
 & NP(A) \ar[ur] & \\
 & \tau_{\geq 0}P(NA) \ar@/_1pc/[uur]_{\tau_{\geq 0}(P)} & 
}
$$ 
of the diagonal dg functor. By lemma~\ref{path} $\tau_{\geq 0}P(NA)$ is a path object of $NA$ in $\dgcatp$. We will show in proposition~\ref{important} that $NP(A)$ is also a path object of $NA$.

\begin{lemma}\label{compar}
Let $A,B\in s\mathbf{Mod}$. The shuffle map $\nabla$ induces a natural surjective chain homotopy equivalence 
$$ N(\underline{s\mathbf{Mod}}(A,B)) \stackrel{\nabla^{\sharp}}{\longrightarrow} \underline{\mathsf{Ch}_{\geq 0}} (NA, NB),$$
which has a natural section induced by the Alexander-Whitney map.
\end{lemma}
\begin{proof}
First note that if $(L,R)$ and $(L',R')$ are adjoint pairs of functors, a natural
transformation $\zeta \colon L \to L'$ induces a natural transformation $\zeta^\sharp \colon R' \to R$
which is a natural equivalence iff $\zeta$ is also.

Fixing a chain complex $NA \in \mathsf{Ch}_{\geq 0}$ let $L,L' \colon 
\mathsf{Ch}_{\geq 0} \to \mathsf{Ch}_{\geq 0}$ be defined by 
\[ L(C):= C \otimes NA, \quad \quad L'(C):= N(\Gamma C \wedge A). \]
Using the Dold-Kan equivalence in the case of $L'$, we see that these functors have right adjoints
\[ R(C) = \underline{\mathsf{Ch}_{\geq 0}} (NA, C), \quad \quad R'(C) = N(\underline{s\mathbf{Mod}}(A,\Gamma C)) \]
respectively. 

The shuffle map determines a natural inclusion $\nabla \colon L \to L'$ which has 
a right inverse given by the Alexander-Whitney map $AW$, see \cite[2.7]{SS}. It follows
that $\nabla^\sharp \colon R' \to R$ is a natural surjection with a section given by $AW^\sharp$.

The fact that $\nabla^\sharp$ is a natural transformation of bi-functors is clear.

Since $\nabla$ is a chain homotopy equivalence, in order to finish the proof it
is now enough to show that the functors $L,L',R,R'$ send chain homotopic maps to chain homotopic maps
(for $(L,R)$ and $(L',R')$ will then induce adjunctions on the homotopy category $\mathsf{Ho}(\mathsf{Ch}_{\geq 0})$ 
and $\nabla \colon L \to L'$ will be a natural isomorphism between endo-functors of $\mathsf{Ho}(\mathsf{Ch}_{\geq 0})$). 

The functors $L$ and $R$ clearly preserve the chain homotopy relation. For the same reason, $L'$ and $R'$ 
preserve the relation on $\mathsf{Ch}_{\geq 0}(C,D)$ defined by the cylinder object 
\[
\xymatrix{
C  \ar[r] \ar[rd] & N(\Gamma C \wedge k\Delta[1]) \ar[d] & \ar[l] \ar[ld] C \\
& C\,. & }
\]
Since the Alexander-Whitney and shuffle maps give maps between this cylinder object and the usual
one, we see that this relation is the usual chain homotopy relation. This concludes the proof.
\end{proof}

We now define a map $\phi$ relating the $\mathsf{Ch}_{\geq 0}$-graphs associated with the dg categories $NP(A)$ and $\tau_{\geq 0}P(NA)$. Observe that:
\begin{itemize}
\item[-] By lemma~\ref{path}, $NP(A)$ and $\tau_{\geq 0}P(NA)$ have exactly the same objects and
\item[-] For each pair of objects $x \stackrel{f}{\rightarrow} y$, $x' \stackrel{f'}{\rightarrow} y'$ in $NP(A)$, the map of lemma~\ref{compar} (with $A=k\Delta[1]$) induces a surjective quasi-isomorphism $\phi_{f,f'}$
$$
\xymatrix{
NA(x,x') \underset{NA(x,y')}{\times} N\underline{s\mathbf{Mod}}(k\Delta[1],A(x,y')) \underset{NA(x,y')}{\times} NA(y,y') \ar@{->>}[d]^-{\mathbf{1} \times \nabla^{\sharp} \times \mathbf{1}}_{\sim} \\
NA(x,x') \underset{NA(x,y')}{\times} \underline{\mathsf{Ch}_{\geq 0}}(Nk\Delta[1],NA(x,y')) \underset{NA(x,y')}{\times} NA(y,y')
}
$$
in $\mathsf{Ch}_{\geq 0}$.
\end{itemize}
\begin{notation}
We denote by
$$ 
\xymatrix{
\phi: NP(A) \ar@{-->}[r] & \tau_{\geq 0}P(NA)
}$$
the map of $\mathsf{Ch}_{\geq 0}$-graphs which is the identity on objects and $\phi_{f,f'}$ on the complexes of morphisms.
\end{notation}
\begin{remark}
Notice that by definition of $P(A)$ and remark~\ref{nova} the map $\phi$ preserve identities and the composition of degree zero morphisms.
\end{remark}

We now establish a `homotopy equivalence lifting property' of $\phi$.
\begin{proposition}\label{lifth}
Let $\alpha$ be a degree zero morphism in $\tau_{\geq 0}P(NA)$ that becomes invertible in $\mathsf{H}_0(\tau_{\geq 0}P(NA))$. Then there exists a degree zero morphism $\overline{\alpha}$ in $NP(A)$ which becomes invertible in $\mathsf{H}_0(NP(A))$ and $\phi(\overline{\alpha})=\alpha$.
\end{proposition}
\begin{proof}
Let 
$$
\alpha: (x \stackrel{f}{\rightarrow} y) \longrightarrow (x' \stackrel{f'}{\rightarrow} y')$$
be a degree zero morphism in $\tau_{\geq 0}P(NA)$. Notice that $\alpha$ is of the form $(m_x, h, m_y)$ with $m_x: x \rightarrow x'$ and $m_y:y \rightarrow y'$ degree zero morphisms in $NA$ and $h: x \rightarrow y$ a degree $1$ morphism in $NA$. Now, by definition of $P(A)$ we can choose a representative $\overline{h} \in A(x,y)_1$ of $h$ and so we obtain a degree zero morphism $\overline{\alpha}=(m_x, \overline{h}, m_y)$ in $NP(A)$ such that
$$
\begin{array}{rcl}
\phi_{f,f'}: NP(A)(f,f') &  \rightarrow & \tau_{\geq 0}P(NA)(f,f')\\
\overline{\alpha}=(m_x, \overline{h}, m_y) & \mapsto & (m_x, h, m_y)\,.
\end{array}
$$
Now suppose that $\alpha$ is invertible in $\mathsf{H}_0(\tau_{\geq 0}P(NA))$. Then there exist morphisms $\beta$ of degree $0$ and $r_1$ and $r_2$ of degree $1$ such that $d(r_1)=\beta \alpha-\mathbf{1}$ and $d(r_2)=\alpha \beta-\mathbf{1}$. As above, we can lift $\beta$ to a morphism $\overline{\beta}$ in $NP(A)$. Since the map $\phi$ preserve the identities and the composition of degree zero morphisms it maps $\overline{\alpha}\overline{\beta}$ to $\alpha \beta$ and $\overline{\beta}\overline{\alpha}$ to $\beta \alpha$. Finally since the maps $\phi_{f,f'}$ are surjective quasi-isomorphisms we can lift $r_1$ to $\overline{r_1}$, resp. $r_2$ to $\overline{r_2}$, in $NP(A)$ by applying the lemma
\cite[2.3.5]{Hovey} to the couple $(r_1,\mathbf{1})$, resp. $(r_2,\mathbf{1})$. This implies that $\overline{\alpha}$ is also invertible in $\mathsf{H}_0(NP(A))$.
\end{proof}

\begin{proposition}\label{important}
In the following commutative diagram in $\dgcatp$
$$
\xymatrix{
NA \ar[rr]^{\Delta} \ar[dr]_{N(I_A)} & & NA\times NA \\
& NP(A) \ar[ur]_{N(P_0)\times N(P_1)}\,, &
}
$$
obtained by applying the normalization functor $N$ to the diagram~(\ref{numero}) in $s\mathsf{Mod}\text{-}\mbox{Cat}$, the dg functor $N(I_A)$ is a quasi-equivalence and $N(P_0)\times N(P_1)$ is a fibration.
\end{proposition}
\begin{proof}
We first prove that $N(I_A)$ is a quasi-equivalence.  By definition of $P(A)$ the dg functor $I_A$ clearly satisfies condition $(1)$ of proposition~\ref{eqfracas}. We now prove that $N(I_A)$ satisfies condition $(ii)$ of definition~\ref{quasieq}. Let $f$ be an object in $NP(A)$. The dg categories $NP(A)$ and $\tau_{\geq 0}P(NA)$ have the same objects and so we can consider $f$ as an object in $\tau_{\geq 0}P(NA)$. Since the dg functor
$$\tau_{\geq 0}(I): NA \longrightarrow \tau_{\geq 0}P(NA)$$
is a quasi-equivalence, see lemma~\ref{path}, there exists an object $x$ in $NA$ and a homotopy equivalence $\alpha$ in $\tau_{\geq 0}P(NA)$ between $I(x)$ and $f$. By proposition~\ref{lifth} we can lift $\alpha$ to a homotopy equivalence $\overline{\alpha}$ in $NP(A)$ and so the dg functor
$$ N(I_A): NA \longrightarrow NP(A)$$
satisfies condition $(ii)$ of definition~\ref{quasieq}. This proves that $N(I)$ is a quasi-equivalence.

We now prove that $N(P_0)\times N(P_1)$ is a fibration. By definition of $P(A)$ the dg functor $N(P_0) \times N(P_1)$ is clearly surjective on the complexes of morphisms. We now prove that it has the right lifting property with respect to the generating trivial cofibration $F$, see definition~\ref{legos}. Let $x \stackrel{f}{\rightarrow} y$ be an object in $NP(A)$ and $\gamma:(x,y) \rightarrow (x',y')$ a homotopy equivalence in $NA \times NA$. Since the dg functor
$$ 
\xymatrix{\tau_{\geq 0}(P): \tau_{\geq 0} P(NA) \ar@{->>}[r] & NA \times NA}$$
is a fibration there exists a homotopy equivalence $\alpha: f \rightarrow f'$ in $\tau_{\geq 0}P(NA)$ such that $\tau_{\geq 0}(P)(\alpha)=\gamma$. Now, by proposition~\ref{lifth} we can lift $\alpha$ to a homotopy equivalence $\overline{\alpha}: f \rightarrow f'$ in $N(PA)$ such that $N(P_0) \times N(P_1) (\overline{\alpha})=\gamma$.

This proves the proposition.
\end{proof}
Notice that the previous proposition implies theorem~\ref{newQuillen}.

\begin{remark}\label{fibrant1}
Since every object in $\dgcatp$ is fibrant, see remark~\ref{fibrant}, all simplicial $k$-linear categories will be fibrant with respect to this Quillen model structure.
\end{remark}

\subsection{Quillen equivalence}
In this subsection we prove that the Quillen adjunction constructed in the previous subsection
$$ \xymatrix{
s\mathbf{Mod}\text{-}\mbox{Cat} \ar@<1ex>[d]^N \\
\dgcatp \ar@<1ex>[u]^L
}
$$
is in fact a Quillen equivalence.

\begin{theorem}\label{equivalenciafinal}
The Quillen adjunction $(L,N)$ is a Quillen equivalence.
\end{theorem}

\begin{proof}
Let $\ca \in \dgcatp$ be a cofibrant dg category and $B$ a simplicial $k$-linear category. Recall from remark~\ref{fibrant1} that every object in $s\mathbf{Mod}\text{-}\mbox{Cat}$ is fibrant. We need to show that a simplicial $k$-linear functor
$$ F: L(\ca) \longrightarrow B$$
is a weak equivalence in $s\mathbf{Mod}\text{-}\mbox{Cat}$ iff the corresponding dg functor
$$ F^{\sharp}: \ca \longrightarrow NB$$
is a quasi-equivalence in $\dgcatp$.

We have the folowing commutative diagram
$$
\xymatrix{
\ca \ar[r]^{F^{\sharp}} \ar[d]_{\eta} & NB \\
NL(\ca) \ar[ur]_{NF} &\,, 
}
$$ 
where $\eta$ is the counit of the adjunction $(L,N)$. Since, by definition, $F$ is a weak equivalence in $s\mathbf{Mod}\text{-}\mbox{Cat}$ iff $NF$ is a quasi-equivalence it is enough to show that $\eta$ is a quasi-equivalence.
The dg functor $\eta$ is the identity map on objects and so it is enough to show that 
$$\eta(x,y): \ca(x,y) \longrightarrow NL(\ca)(x,y),\,\,\,\,x,y \in \ca$$
is a quasi-isomorphism. 
Now, let $I$ be the set of objects of $\ca$. Since $\ca$ is cofibrant in $\dgcatp$ it clearly stays cofibrant when considered as an object of the Quillen model structure on $\mathsf{Ch}_{\geq 0}^I\text{-}\mbox{Cat}$, see remark~\ref{Qmodel}. By proposition $6.4$ of \cite{SS} the adjunction morphism in $\mathsf{Ch}_{\geq 0}^I\text{-}\mbox{Gr}$
$$ \Gamma U (\ca) \longrightarrow L_I(\ca)$$
is such that
$$ \Gamma U(\ca)(x,y) \longrightarrow L_I(\ca)(x,y)$$
induces an isomorphism in $\pi_i$ for $i\geq 0$ and for all objects $x,y \in \Gamma U (\ca)$. This implies by the Dold-Kan equivalence that
$$\ca(x,y)= N(\Gamma U(\ca)(x,y)) \stackrel{\sim}{\longrightarrow} N(L_I(\ca)(x,y)), \,\,\,\,x,y \in \ca$$
is a quasi-isomorphism and so 
$$\eta(x,y): \ca(x,y) \longrightarrow NL\ca(x,y), \,\,\,\, x, y \in \ca$$ is a quasi-isomorphism.
This proves the theorem. 
\end{proof}

\begin{remark}
Notice that the objects in $\dgcatp$, resp. in $s\mathbf{Mod}\text{-}\mbox{Cat}$, with only one object consist exactly on the connective dg algebras, see \cite[1.1]{SS}, resp. simplicial $k$-algebras. We have the following commutative diagram
$$
\xymatrix{
*+<1pc>{s\mathbf{Alg}} \ar@<1ex>[d]^N \ar@{^{(}->}[r] & s\mathbf{Mod}\text{-}\mbox{Cat} \ar@<1ex>[d]^N \\
\mathbf{DGA}_{\geq 0} \ar@<1ex>[u]^L \ar@{^{(}->}[r] & \dgcatp \ar@<1ex>[u]^L\,,
}
$$
where $\mathbf{DGA}_{\geq 0}$ denotes the category of connective dg algebras and $s\mathbf{Alg}$ the category of simplicial $k$-algebras. Observe that if we restrict the Quillen model structures to these full subcategories we obtain Shipley-Schwede's Quillen equivalence \cite[1.1]{SS}
$$
\xymatrix{
s\mathbf{Alg} \ar@<1ex>[d]^N \\
\mathbf{DGA}_{\geq 0} \ar@<1ex>[u]^L\,.
}
$$
\end{remark}
We have then extended Shipley-Schwede's work to a `several objects' context: the notion of weak equivalence in $s\mathbf{Mod}\text{-}\mbox{Cat}$ and $\dgcatp$ (see definition~\ref{quasieq} and proposition~\ref{abmorph}) is now a mixture between quasi-isomorphisms and categorical equivalences. 

\section{The global picture}
Recall from \cite[III]{Jardine} that we have an adjunction
$$
\xymatrix{
s\mathbf{Mod} \ar@<1ex>[d]^U \\
s\mathbf{Set} \ar@<1ex>[u]^{k(-)},
}
$$
where $U$ is the forgetful functor and $k(-)$ the $k$-linearization functor. The functor $k(-)$ is lax strong monoidal and so we have the natural adjunction
$$
\xymatrix{
s\mathbf{Mod}\text{-}\mbox{Cat} \ar@<1ex>[d]^U \\
s\mathbf{Set}\text{-}\mbox{Cat} \ar@<1ex>[u]^{k(-)}.
}
$$
Recall from \cite[1.1]{Bergner} that the category $s\mathbf{Set}\text{-}\mbox{Cat}$ is endowed with a Quillen model structure whose weak equivalences are the Dwyer-Kan (=DK) equivalences. Let us recall this notion.

\begin{definition}
A simplicial functor $F:A \rightarrow B$ is a {\em Dwyer-Kan equivalence} if:
\begin{itemize}
\item[-] for any objects $x$ and $y$ in $A$, the map
$$ F(x,y):A(x,y)\longrightarrow B(Fx,Fy)$$
is a weak equivalence of simplicial sets and
\item[-] the induced functor $\pi_0(F): \pi_0(A) \longrightarrow \pi_0(B)$ is essentially surjective.
\end{itemize}
\end{definition}

\begin{proposition}
The adjunction $(k(-), U)$ is a Quillen adjunction, when we consider on $s\mathbf{Mod}\text{-}\mbox{Cat}$ the Quillen model structure of theorem~\ref{newQuillen}.
\end{proposition}

\begin{proof}
We first observe that by proposition~\ref{eqfracas} the functor
$$U : s\mathbf{Mod}\text{-}\mbox{Cat} \longrightarrow s\mathbf{Set}\text{-}\mbox{Cat}$$
preserves weak equivalences.

We now show that it also preserves fibrations. Let $G:A \rightarrow B$ be a simplicial $k$-linear functor such that $NG: NA \rightarrow NB$ is a fibration in $\dgcatp$. We need to show that $UG$ is a fibration in $s\mathbf{Set}\text{-}\mbox{Cat}$. Recall from \cite{Bergner} that $UG$ is a fibration iff:
\begin{itemize}
\item[(F1)] for any object $x$ and $y$ in $UA$, the map
$$UG(x,y): UA(x,y) \longrightarrow UB(Gx,Gy)$$
 is a fibration in $s\mathbf{Set}$ and 
\item[(F2)] for any object $x\in UA$, $y \in UB$ and homotopy equivalence $f: Gx \rightarrow y$ in $UB$ (= $f$ becomes invertible in $\pi_0(UB)$), there is an object $z \in A$ and a homotopy equivalence $h:x \rightarrow z$ in $UA$ such that $UG(h)=f$. 
\end{itemize}
Since by hypothesis $NG:NA \rightarrow NB$ is a fibration in $\dgcatp$, the dg functors $R(n), n\geq 1$ (which belong to the set $J$ of generating trivial cofibrations) allow us to conclude that the morphisms
$$ NG(x,y)_n: NA(x,y)_n \longrightarrow NB(Gx,Gy)_n, \,\,\,\, x,y \in A$$
are surjective for $n \geq 1$. Now, by \cite[III-2.11]{Jardine}, $UG(x,y)$ is a fibration in $s\mathbf{Set}$ iff the morphisms $NG(x,y)_n$ are surjective for $n \geq 1$. This implies that condition $(F1)$ is verified.

Concerning condition (F2), let $x \in UA$, $y \in UB$ and $f: Gx \rightarrow y$ be a homotopy equivalence in $UB$. This means that $f$ is invertible in $\pi_0(B)$ and so by lemma~\ref{key} $N(f)$ is also invertible in $\mathsf{H}_0(NB)$. This data allow us to construct, using proposition $1.7$ in \cite{Tab}, the following (solid) commutative square
$$
\xymatrix{
*+<1pc>{\underline{k}} \ar[r] \ar@{>->}[d]_F^{\sim} & NA \ar@{->>}[d]^{NG}\\
\mathcal{K} \ar[r] \ar@{.>}[ur] & NB\,.
}
$$
Since $NG$ is a fibration in $\dgcatp$ we can lift $Nf$ to a morphism $h:x \rightarrow z$ in $NA$ which is invertible in $H_0(NA)$. Since the $0$-simplex morphisms in $A$ and the degree zero morphisms in $NA$ are exactly the same, lemma~\ref{key} implies that $h: x \rightarrow z$, when considered as a morphism in $UA$, satisfies condition (F2).

This proves the proposition.
\end{proof}

We have obtained the following zig-zag of Quillen adjunctions relating the homotopy theories of differential graded and simplicial categories:
$$
\xymatrix{
s\mathbf{Set}\text{-}\mbox{Cat} \ar@<-1ex>[d]_{k(-)} \\
s\mathbf{Mod}\text{-}\mbox{Cat} \ar@<-1ex>[u]_U \ar@<1ex>[d]^N \\
\dgcatp \ar@<1ex>[u]^L  \ar@<-1ex>@{_{(}->}[d]\\
\dgcat \ar@<-1ex>[u]_{\tau_{\geq 0}}.
}
$$
\begin{remark}\label{Simprk}
Since the adjunction $(L,N)$ is a Quillen equivalence
\begin{itemize}
\item[-] the composed functor $\mathbb{L}(N\circ k(-)): \mathsf{Ho}(s\mathbf{Set}\text{-}\mbox{Cat}) \longrightarrow \mathsf{Ho}(\dgcat)$ preserves homotopy colimits and
\item[-] the composed functor $\mathbb{R}(U\circ L\circ \tau_{\geq 0}): \mathsf{Ho}(\dgcat) \longrightarrow \mathsf{Ho}(s\mathbf{Set}\text{-}\mbox{Cat})$ preserves homotopy limits. 
\end{itemize}
\end{remark}
The following result was proved by Simpson in an adhoc way in \cite[5.1]{Simpson}.
\begin{corollary}
Let $F:\ca \rightarrow \cb$ be a dg functor and $b$ an object of $\cb$. Then the homotopy fiber of $F$ over $b$, denoted by $\mathbf{HFib}(F)/b$, is equivalent to the homotopy fiber of $\mathbb{R}(U\circ L\circ \tau_{\geq 0})(F)$ over $b$:
$$ \mathbf{HFib}(F)/b \stackrel{\sim}{\longrightarrow} \mathbf{HFib}(NF)/b \,.$$
\end{corollary}
\begin{proof}
It follows from remark~\ref{Simprk}.
\end{proof}

\end{document}